\documentclass[a4paper,11pt]{article}
\usepackage[mathscr]{eucal}
\usepackage{amssymb}
\usepackage{latexsym}
\usepackage{theorem}
\theoremstyle{plain}

\makeatletter
  
   \@addtoreset{equation}{section}
\makeatother     
\renewcommand{\baselinestretch}{1}

\def\<{\left<} \def\>{\right>}
\newtheorem{theorem}{Theorem}
\newtheorem{lemma}[theorem]{Lemma}

\newtheorem{remark}[theorem]{Remark}
\newtheorem{defn}[theorem]{Definition}
\def\proof{\noindent{\it Proof: }}
\def\bea{\begin{eqnarray} }
\def\eea{\end{eqnarray} }
\def\be{\begin{equation} }
\def\ee{\end{equation} }
\def\qed{\ifhmode\unskip\nobreak\fi\ifmmode\ifinner\else\hskip5pt
\fi\fi\hbox{\hskip5 pt \vrule width4 pt height6 pt depth1.5 pt \hskip1pt }}

\begin{document}

\title{Quasi-biharmonic 
Lagrangian surfaces in Lorentzian complex space forms 
\footnote{Annali di Matematica Pura ed Applicata {\bf 192} (2013), 191-201.
A section is added at the end.}}
\author{Toru Sasahara}

\date{}
\maketitle

\begin{abstract}
{\footnotesize 
In this paper, we introduce the notion of a quasi-biharmonic submanifold in a pseudo-Riemannian manifold
and classify quasi-biharmonic marginally trapped Lagrangian surfaces in Lorentzian complex space
forms.}

\end{abstract}




\section{Introduction}

A submanifold with lightlike mean curvature vector field  is  called a {\it marginally trapped} submanifold.
In the theory of cosmic black holes, a marginally trapped surface in a space-time plays an
extremely important role.  Recently, some classification results on marginally trapped surfaces from
the viewpoint of differential geometry have been obtained (see, for instance, \cite{ch3}).

On the other hand, a submanifold 
 is called {\it biharmonic} if the bitension field of the isometric immersion defining the submanifold vanishes identically.
The theory of biharmonic submanifolds
has advanced greatly during this last decade 
(see, for instance,  \cite{bal}  and \cite{MO}).  
 This paper introduces the notion of a {\it quasi-biharmonic} submanifold.
It is a submanifold such that
the bitension field of the isometric immersion defining the submanifold is lightlike at each point.

In \cite{sasa2}, the author has classified biharmonic marginally trapped Lagrangian surfaces in
Lorentzian complex space forms. 
They exist only in the flat Lorentzian complex plane. 
In this paper, we classify
quasi-biharmonic marginally trapped Lagrangian surfaces in 
Lorentzian complex space forms. 
We find that the situation in the quasi-biharmonic case
 is quite different from the biharmonic case.
In fact,
 there exist a lot of quasi-biharmonic marginally trapped Lagrangian surfaces
in nonflat Lorentzian complex space forms.

\section{Preliminaries}

\subsection{Lagrangian submanifolds in complex space forms} 

Let $\tilde M^n_s(4\epsilon)$ be a complex space form of complex dimension $n$,
complex index $s(\geq 0)$ and constant holomorphic sectional
curvature $4\epsilon$. The complex index is defined as the complex dimension of
the largest complex negative definite vector subspace of the tangent space.
If $s=1$, it is called {\it Lorentzian}.
The curvature tensor $\tilde R$ of $\tilde M^n_s(4\epsilon)$ is given by
\bea
\tilde R(X, Y)Z&=&\epsilon\{\<Y, Z\>X-\<X, Z\>Y+\<Z, JY\>JX\nonumber\\
&&-\<Z, JX\>JY+2\<X, JY\>JZ\},\eea
where $\< , \>$ and $J$ are the metric tensor and
the almost complex structure of $\tilde M^n_s(4\epsilon)$ respectively.

Let ${\bf C}^n_s$ be the $n$-dimensional complex space with complex coordinates $z_1, \ldots, z_n$, endowed with
the metric
$g_{n, s}(z, w)={\rm Re}(-\sum_{j=1}^sz_j\bar{w}_j+\sum_{i=s+1}^nz_i\bar{w}_i)$.
Put $S^{2n+1}_{2s}(\epsilon)=\{z\in {\bf C}^n_s: g_{n+1, s}(z, z)=\frac{1}{\epsilon}\}$ for $\epsilon>0$
and $H^{2n+1}_{2s+1}(\epsilon)=\{z\in {\bf C}^n_s : g_{n+1, s+1}(z, z)=\frac{1}{\epsilon}\}$ for $\epsilon<0$.

The Hopf fibrations 
\bea
\pi: S^{2n+1}_{2s}(\epsilon)\rightarrow {\bf C}P^n_s(4\epsilon): z\rightarrow z\cdot {\bf C}^{*},\nonumber\\
\pi: H^{2n+1}_{2s+1}(\epsilon)\rightarrow {\bf C}H^n_s(4\epsilon): z\rightarrow z\cdot {\bf C}^{*},\nonumber
\eea
give ${\bf C}P^n_s(4\epsilon)$ and ${\bf C}H^n_s(4\epsilon)$  a unique pseudo-Riemannian metric of complex index $s$  and curvature tensor (2.1) 
such that
$\pi$ is a pseudo-Riemannian submersion  respectively.

Barros and Romero \cite{ba 2} showed that locally any  complex space form
$\tilde M^n_s(4\epsilon)$
is isometric holomorphically to ${\bf C}^n_s$, ${\bf C}P^n_s(4\epsilon)$ ${\bf C}H^n_s(4\epsilon)$
according to $\epsilon=0$, $\epsilon>0$ or $\epsilon<0$. 

An $n$-dimensional
submanifold $M$
isometrically immersed in $\tilde M^n_s(4\epsilon)$ is called {\it Lagrangian} if $J$ interchanges the tangent  and the normal spaces of $M$.
For a Lagrangian submanifold  $M$ of 
complex space form $\tilde M^n_s(4\epsilon)$,
we denote by $\nabla$ and
$\tilde\nabla$ the Levi-Civita connections on $M$ and 
$\tilde M^n_s(4\epsilon)$, respectively. The
formulas of Gauss and Weingarten are given respectively  by
\bea
& &\tilde \nabla_XY= \nabla_XY+h(X,Y),\\
& &\tilde\nabla_X JY= -A_{JY}X+D_XJY,\eea
for $X, Y$ tangent to $M$, where $h,A$ and $D$ are the second fundamental
form, the shape operator and the normal
connection.
The mean curvature vector field is defined by 
$H=\frac{1}{n}{\rm tr}h$. 
The shape operator and the second fundamental form are related by
\bea
\<h(X, Y), JZ\>=\<A_{JZ}X, Y\>
\eea
for $X, Y, Z$ tangent to $M$.
Since $J$ is parallel, by (2.2) and (2.3) we have
\bea
& &D_XJY=J(\nabla_XY),\\
& &A_{JY}X=-Jh(X,Y)=A_{JX}Y.
\eea

The equations of
Gauss, Codazzi are given respectively by
\be
\<R(X,Y)Z,W\>=
\epsilon(\<X,W\>\<Y,Z\>-\<X,Z\>\<Y,W\>)+\<[A_{JZ},
A_{JW}](X),Y\>,
\ee
\be
({\bar\nabla}_{X}h)(Y,Z)=
({\bar\nabla}_{Y}h)(X,Z),
\ee
where $X,Y,Z,W$ 
 are vectors tangent  to
$M$, 
$R(X, Y)=[\nabla_X, \nabla_Y]-\nabla_{[X, Y]}$ and 
$\bar\nabla h$ is defined by
\be ({\bar\nabla}_{X}h)(Y,Z)= D_X h(Y,Z) - h(\nabla_X
Y,Z) - h(Y,\nabla_X Z).
\ee

\subsection{Legendre curves in the light cone}

A vector $X$ is called {\it spacelike} (resp. {\it timelike}) if it satisfies
$\<X, X\>>0$ (resp. $\<X, X\><0$).
A vector $X$ is called {\it lightlike}
if it is nonzero and it satisfies $\<X, X\>=0$.
The light cone $\mathcal{L}C\subset  {\bf C}_1^2$ is defined by $\mathcal{L}C$=$\{z\in{\bf C}_s^n: \<z, z\>=0\}$.
A curve $z(t)$ is called {\it null} if $z^{\prime}$ is lightlike for any $t$.

A 
curve $z(t)$ in $\mathcal{L}$C
is called {\it Legendre} if $\<z^{\prime}, iz\>=0$ holds  for any $t$.
A Legendre curve $z(t)$ in $\mathcal{L}C$ is called {\it special Legendre} if  $\<iz^{\prime}, z^{\prime\prime}\>=0$ holds for any $t$.
For a unit speed special Legendre curve $z(t)$ in $\mathcal{L}C$, the {\it squared Legendre curvature} and the {\it Legendre torsion}
are defined by $\hat\kappa^2=\<z^{\prime}, z^{\prime}\>\<z^{\prime\prime}, z^{\prime\prime}\>$
and $\hat\tau=\<z^{\prime}, z^{\prime}\>\<z^{\prime\prime}, iz^{\prime\prime\prime}\>$, respectively.
For further details on Legendre curves in the light cone, see \cite{ch2} and \cite{chd}.

\subsection{Bitension field}

Let $M$ and  $N$ be pseudo-Riemannian manifolds of dimension $m$ and $n$ respectively,
and $\phi:M \to N$ a smooth map.
We denote by $\nabla$ and $\tilde\nabla$ the Levi-Civita connections on $M$ and $N$ respectively.
Then the {\it tension field}
$\tau(\phi)$ is a section
of the vector bundle $\phi^{*}TN$
defined by
$$
\tau(\phi)=\mathrm{tr}(\nabla^{\phi} d \phi)=
\sum_{i=1}^{m}\<e_i, e_i\>\{\nabla^{\phi}_{e_i}d\phi(e_i)
-d\phi(\nabla_{e_i}e_i)\}.
$$
Here $\nabla^{\phi}$ and $\{e_i\}$ denote 
the induced connection by $\phi$ on the bundle $\phi^{*}TN$,  which is the
pull-back of $\tilde\nabla$,
and a local orthonormal
frame of $M$,
respectively.
If $\phi$ is an isometric immersion, then $\tau(\phi)$ and the mean curvature vector field $H$ of $M$ are related by
\be
\tau(\phi)=mH.
\ee

A smooth map $\phi$ is said to be 
a {\it harmonic map} if $\tau(\phi)=0$ at each point on $M$.
If $M$ and $N$ are Riemannian manifolds, then $\phi$ is harmonic if and only if
it is a critical 
point of the {\it energy}
$$
E(\phi)=\int_{\Omega}|d\phi|^2dv_{g}
$$
over every compactly supported region $\Omega$ of $M$.

We define the {\it bitension field} as 
\be
\tau_2(\phi)=\sum_{i=1}^{m}\<e_i, e_i\>\{
(\nabla^{\phi}_{e_i}\nabla^{\phi}_{e_i}-\nabla^{\phi}_{\nabla_{e_i}e_i})\tau
+R^{N}(\tau,d\phi (e_i))d\phi(e_i)\},
\ee
where $R^N$ is the curvature tensor of $N$.
If $\phi$ is an isometric immersion and $N$ is the complex space form $\tilde M^n_s(4\epsilon)$,
 then (2.1), (2.10) and (2.11) yield
\bea
\tau_2(\phi)=-m\Delta H+5m\epsilon H,
\eea
where $\Delta=-\sum_{i=1}^{m}\<e_i, e_i\>(\tilde\nabla_{e_i}\tilde\nabla_{e_i}-\tilde\nabla_{\nabla_{e_i}e_i})$.

A smooth map $\phi$ is called {\it biharmonic} if $\tau_2(\phi)=0$ at each point on $M$.
Harmonic maps are clearly biharmonic. 
When $M$ and $N$ are Riemannian manifolds, the biharmonic map $\phi$ is characterized as 
a critical point of the {\it bienergy}
$$
E_{2}(\phi)=\int_{\Omega}|\tau(\phi)|^2dv_{g},
$$
over every compactly supported region $\Omega$ of $M$.
For recent information on biharmonic maps between Riemannian manifolds, see \cite{bal} and \cite{MO}.


A pseudo-Riemannian submanifold 
in a pseudo-Riemannian manifold isometrically immersed by $\phi$ 
is called {\it marginally trapped} (or {\it quasi-minimal}) if the mean curvature vector field is lightlike, equivalently,
 $\tau(\phi)$ is lightlike at each point on the submanifold.
Analogously, we introduce the new class of submanifolds in  pseudo-Riemannian manifolds as follows. 
\begin{defn}
{\rm A pseudo-Riemannian submanifold in a pseudo-Riemannian manifold isometrically immersed by $\phi$ 
is called} {\it quasi-biharmonic} {\rm if $\tau_2(\phi)$ is lightlike at each point on the submanifold.}
\end{defn}

\section{The bitension field of marginally trapped Lagrangian immersions}

Let $\phi: M\rightarrow \tilde M^2_1(4\epsilon)$ be a Lagrangian isometric immersion into a
2-dimensional Lorentzian complex space forms. Then, the real index of $M$ is one.
Choose a local orthonormal frame $\{e_1, e_2\}$
such that  $\<e_1, e_1\>=1$ and $\<e_2, e_2\>=-1$.
Put $\omega_i^j(e_k)=\langle \nabla_{e_k}e_i,e_j\rangle\langle {e_j},e_j\rangle$ for $i, j, k=1,2$.
Note that
 \be
\omega_1^2(e_k)=\omega_2^1(e_k), \quad \omega_1^1(e_k)=\omega_2^2(e_k)=0.
\ee
 
It follows from (2.4) and (2.6) that the
 second fundamental form and the shape operator take the form
\bea
&&h(e_1, e_1)=aJe_1+bJe_2,\\
&&h(e_1, e_2)=-bJe_1+cJe_2,\\ 
&&h(e_2, e_2)=-cJe_1+dJe_2,\\
&&A_{Je_1}e_1=ae_1+be_2,\\
&&A_{Je_1}e_2=-be_1+ce_2,\\ 
&&A_{Je_2}e_1=-be_1+ce_2,\\
&&A_{Je_2}e_2=-ce_1+de_2,
\eea
 for some functions $a$, $b$, $c$, $d$.

We compute (2.9) using (3.2)-(3.4).
In view of (3.1)
we get 
\bea
(\bar\nabla_{e_2}h)(e_1,e_1)&=&(e_2a+3b\omega_1^2(e_2))Je_1\nonumber\\
& &+(e_2b+(a-2c)\omega_1^2(e_2))Je_2, \\
(\bar\nabla_{e_1}h)(e_1,e_2)
&=&-(e_1b+(a-2c)\omega_1^2(e_1))Je_1\nonumber\\
& &+(e_1c-(2b+d)\omega_1^2(e_1))Je_2,\\
(\bar\nabla_{e_1}h)(e_2,e_2)&=&(-e_1c+(2b+d)\omega_1^2(e_1))Je_1\nonumber\\
& &+(e_1d-3c\omega_1^2(e_1))Je_2,\\
(\bar\nabla_{e_2}h)(e_1,e_2)&=&-(e_2b+(a-2c)\omega_1^2(e_2))Je_1\nonumber\\
& &+(e_2c-(2b+d)\omega_1^2(e_2))Je_2.
\eea
From (2.8) and (3.9)-(3.12) we obtain 
\bea
& &e_2a+3b\omega_1^2(e_2)=-e_1b-(a-2c)\omega_1^2(e_1),\\
& &e_2b+(a-2c)\omega_1^2(e_2)=e_1c-(2b+d)\omega_1^2(e_1),\\
& &e_1d-3c\omega_1^2(e_1)=e_2c-(2b+d)\omega_1^2(e_2).
\eea

Denote the Gauss curvature of $M$ by $G$. Then, 
the Gauss equation (2.7) implies that $G$ satisfies
\be
G=ac+b^2+bd-c^2+\epsilon.
\ee

From now on, 
let us assume that $\phi: M\rightarrow \tilde M_1^2(4\epsilon)$ is a marginally trapped Lagrangian immersion 
into a
2-dimensional Lorentzian complex space forms and
$\{e_1, e_2\}$ is an orthonormal frame on $M$ satisfying
$\<e_1, e_1\>=1$, $\<e_2, e_2\>=-1$, 
\be
H=\alpha(Je_1+Je_2)
\ee
for some  function $\alpha$ which is nowhere zero, 
and the second fundamental form is expressed as $(3.2)-(3.4)$.

The main purpose of this section is to
express the bitension of $\phi$ by $a$, $b$, $c$, $\alpha$, $e_1$, $e_2$ and $\omega_i^j(e_k)$.

By (3.2), (3.4) and (3.17) we have 
\be
2\alpha=a+c=b-d.
\ee
Combining (3.16) and (3.18) yields
\bea
G=(a-2b-c)(c-b)+\epsilon.
\eea

In view of (2.12), we need the following basic formula, which is obtained by a straightforward computation 
using (2.2) and (2.3) (cf. (4.7) in p.270 of \cite{ch1}):
\be
\Delta H=\Delta^DH+\sum_{i=1}^2\<e_i, e_i\>h(e_i, A_He_i)+\sum_{i=1}^2\<e_i, e_i\>(A_{D_{e_i}H}e_i+(\nabla_{e_i}
A_H)e_i), 
\ee
where
$\Delta^D=-\sum_{i=1}^{2}
\<e_i, e_i\>(D_{e_i}D_{e_i}-D_{\nabla_{e_i}e_i})$.

 First, we  have
\begin{lemma} 
Let $M$ be a marginally trapped Lagrangian surface in $\tilde M_1^2(4\epsilon)$. Then, the normal part of $\Delta H$ is expressed as
\bea
(\Delta H)^{\perp}&=&\alpha\{(a-2b-c)(a+2b-c)-\epsilon\}Je_1\nonumber\\
&&+\alpha\{(a-2b-c)(-a+2b-3c)-\epsilon\}Je_2.
\eea
\end{lemma}
\proof By  Proposition 3 in \cite{sasa1}, we have 
\be
\Delta^{D}H=-GH.
\ee
( In \cite{sasa1},  the relation $\Delta^{D}H=GH$ was derived. But the sign of the Gauss curvature (3.16) in \cite{sasa1} was incorrect.)
Using (3.2)-(3.4) we get 
\be
\sum_{i=1}^2\<e_i, e_i\>h(e_i, A_He_i)=\alpha(a-2b-c)\{(a+b)Je_1+(-a+b-2c)Je_2\}.
\ee
From (3.19), (3.20) (3.22) and (3.23) we obtain (3.21).
\qed\vspace{1.5ex}

Next we compute the tangential part of $\Delta H$.
To do so, we need the following Lemma.
\begin{lemma} Let $M$ be a marginally trapped Lagrangian surface in $\tilde M_1^2(4\epsilon)$. Then, we have
\bea 
&&e_1\alpha+e_2\alpha=-\alpha(\omega_1^2(e_1)+\omega^2_1(e_2)),\\
&&e_1a-e_1b+e_2b+e_2c=(a-3b-2c)\omega_1^2(e_1)+(-2a+3b+c)\omega_1^2(e_2),\\
&&e_2a+e_1b-e_2b+e_1c=
(-2a+3b+c)\omega_1^2(e_1)+(a-3b-2c)\omega_1^2(e_2).
\eea
\end{lemma}
\proof
Combining (3.13), (3.15) and (3.18) shows (3.24). 
Similarly, by (3.14), (3.15) and (3.18) we get (3.25).
 We obtain (3.26) from (3.24), (3.25) and (3.18) immediately.
 \qed\vspace{1.5ex}


\begin{lemma}Let $M$ be a marginally trapped Lagrangian surface in $\tilde M_1^2(4\epsilon)$. Then,
the tangential part of $\Delta H$ is expressed as
\bea
(\Delta H)^{\top}=2(a-2b-c)\{e_1\alpha+\alpha\omega_1^2(e_1)\}(e_1-e_2).
\eea
\end{lemma}

\proof
It follows from (3.5)-(3.8) and (3.18) that
\bea
&&A_He_1=\alpha\{(a-b)e_1+(b+c)e_2\},\\
&&A_He_2=\alpha\{-(b+c)e_1+(b-a)e_2\}.
\eea
By (3.1) (3.28) and (3.29) we get
\bea
\nabla_{e_1}(A_He_1)&=&(e_1\alpha)\{(a-b)e_1+(b+c)e_2\}\nonumber\\
&&+\alpha\{(e_1a-e_1b)e_1+(e_1b+e_1c)e_2\nonumber\\
&&+(a-b)\omega_1^2(e_1)e_2
+(b+c)\omega_1^2(e_1)e_1\},\\
A_H(\nabla_{e_1}e_1)&=&\alpha\omega_1^2(e_1)\{-(b+c)e_1+(b-a)e_2\},\\
\nabla_{e_2}(A_He_2)&=&(e_2\alpha)\{-(b+c)e_1+(b-a)e_2\}\nonumber\\
&&+\alpha\{(-e_2b-e_2c)e_1+(e_2b-e_2a)e_2\}\nonumber\\
&&+\alpha\{-(b+c)\omega_1^2(e_2)e_2+(b-a)\omega_1^2(e_2)e_1\},\\
A_H(\nabla_{e_2}e_2)&=&\alpha(\omega_1^2(e_2))\{(a-b)e_1+(b+c)e_2\}.
\eea

Using (2.5), (3.1) and (3.17) we obtain 
\bea
D_{e_1}H=\{e_1\alpha+\alpha\omega_1^2(e_1)\}(Je_1+Je_2),\nonumber\\
D_{e_2}H=\{e_2\alpha+\alpha\omega_1^2(e_2)\}(Je_1+Je_2),\nonumber
\eea
which imply that
\bea
A_{D_{e_1}H}e_1=\{e_1\alpha+\alpha\omega_1^2(e_1)\}\{(a-b)e_1+(b+c)e_2\},\\
A_{D_{e_2}H}e_2=\{e_2\alpha+\alpha\omega_1^2(e_2)\}\{-(b+c)e_1+(b-a)e_2\}.
\eea

By a straightforward computation using  (3.20), (3.30)-(3.35) we find that
the tangential part $(\Delta H)^{\top}$ of $\Delta H$ satisfies
\bea
\<(\Delta H)^{\top},e_1\>&=&(e_1\alpha)(a-b)+(e_2\alpha)(b+c)+\alpha(e_1a-e_1b+e_2b+e_2c)\nonumber\\
&&+2\alpha(b+c)\omega_1^2(e_1)
+2\alpha(a-b)\omega_1^2(e_2)\nonumber\\
&&+\{e_1\alpha+\alpha\omega_1^2(e_1)\}(a-b)+\{e_2\alpha+\alpha\omega_1^2(e_2)\}(b+c),\\
-\<(\Delta H)^{\top},e_2\>&=&(e_1\alpha)(b+c)-(e_2\alpha)(b-a)+\alpha(e_2a+e_1b-e_2b+e_1c)\nonumber\\
&&+2\alpha(a-b)\omega_1^2(e_1)
+2\alpha(b+c)\omega_1^2(e_2)\nonumber\\
&&+\{e_1\alpha+\alpha\omega_1^2(e_1)\}(b+c)-\{e_2\alpha+\alpha\omega_1^2(e_2)\}(b-a).
\eea
Substituting  (3.24)-(3.26) into (3.36) and (3.37) gives us (3.27).
\qed\vspace{1.5ex}

By (2.12), Lemma 2 and Lemma 4,  we obtain the following.
\begin{lemma}Let $\phi: M\rightarrow \tilde M_1^2(4\epsilon)$ be a marginally trapped Lagrangian immersion.
The bitension field of $\phi$ is expressed as
\bea
\tau_2(\phi)&=&-4(a-2b-c)\{e_1\alpha+\alpha\omega_1^2(e_1)\}(e_1-e_2)\nonumber\\
&&-2\alpha\{(a-2b-c)(a+2b-c)-6\epsilon\}Je_1\nonumber\\
&&-2\alpha\{(a-2b-c)(-a+2b-3c)-6\epsilon\}Je_2.
\eea
\end{lemma}






\section{Main results} 
Recently, the following result has been obtained by the author.
\begin{theorem}{\bf(\cite{sasa2})}
Let $M$ be a biharmonic marginally trapped Lagrangian surface in a 
$2$-dimensional Lorentzian complex space form of constant holomorphic sectional curvature $4\epsilon$.
Then $\epsilon=0$, that is, the ambient space is ${\bf C}_1^2$. Moreover, $M$ is locally congruent to 
$$
\phi(x, y)=c_1xe^{if(y)}+z(y),\nonumber
$$
where $f(y)$ is a real-valued function, $c_1$ is a lightlike vector,
$z(y)$ is a null curve in ${\bf C}_1^2$ satisfying $\<iz^{\prime}, c_1e^{if(y)}\>=0$ and
$\<z^{\prime}, c_1e^{if(y)}\>=-1$.
\end{theorem}

In this section, we classify  
{ quasi-biharmonic} marginally trapped Lagrangian surfaces in 2-dimensional Lorentzian complex space forms. Very unlike the biharmonic case,
there exist  a lot of quasi-biharmonic marginally 
trapped Lagrangian surfaces in 2-dimensional {\it nonflat}  Lorentzian complex space forms.

In the case when the ambient space is flat, we have
\begin{theorem}
Let $M$ be a quasi-biharmonic marginally trapped Lagrangian surface in ${\bf C}_1^2$.
Then, $M$  is locally congruent to
\be
\phi (x, y)=e^{i\mu y}z(x), 
\ee
where $\mu$ is a  nonzero real number  and $z(x)$ is a null curve in the light cone $\mathcal{L}$C satisfying 
\bea
\<z, iz^{\prime}\>=\mu^{-1}, \quad z^{\prime\prime}\ne 0
\eea
at each point on $M$.
\end{theorem}
\proof Let $\phi:M\rightarrow{\bf C}_1^2$ be a marginally trapped Lagrangian immersion
and let $\{e_1, e_2\}$ be  an orthonormal frame on $M$ satisfying
$\<e_1, e_1\>=1$, $\<e_2, e_2\>=-1$ and (3.17). Suppose that the second fundamental form is given by (3.2)-(3.3).

If $M$ is quasi-biharmonic, that is, $\tau_2(\phi )$ is lightlike, then from (3.38)
we obtain that $a-2b-c\ne 0$ and $b=c$. 
Hence, we get $G=0$ by (3.19). 
Thus, there exists a local  coordinate system $\{s, t\}$
such that the metric tensor of $M$  is given by
\be
g=ds^2-dt^2.
\ee
We may assume that $e_1=\partial_s$ and $e_2=\partial_t$. 
We put $\partial_x=\frac{1}{\sqrt{2}}(\partial_s-\partial_t)$ and 
$\partial_y=-\frac{1}{\sqrt{2}}(\partial_s+\partial_t)$.
Then,  by (4.3) the metric tensor $g$ is expressed as
\be
g=-dxdy.
\ee
Moreover, a straightforward computation using (3.2)-(3.4) shows 
\bea
&&h(\partial_x, \partial_x)=\frac{a-3b}{\sqrt{2}}J\partial_y,\nonumber\\
&&h(\partial_x, \partial_y)=\sqrt{2}\alpha J\partial _x,\\
&&h(\partial_y, \partial_y)=\sqrt{2}\alpha J\partial _y.\nonumber
\eea

Put $\lambda=\frac{a-3b}{\sqrt{2}}$ and $\mu=\sqrt{2}\alpha$. Note that
$\lambda\ne 0$ and $\mu\ne 0$ at each point of $M$.
By (2.1), (4.4) and (4.5), we see that the Lagrangian immersion $\phi(x, y)$ satisfies
\be
\phi_{xx}=i\lambda \phi_y, \quad \phi_{xy}=i\mu \phi_x, \quad \phi_{yy}=i\mu \phi_y.
\ee
The compatibility condition of (4.6) is given by
\be
\lambda_y=\mu_x=0, \quad \mu_y=0.
\ee
From (4.7) we obtain that $\lambda=\lambda(x)$ and $\mu$ is a nonzero
constant.

By solving the last two equations of (4.6), we have
\be
f(x,y)=e^{i\mu y}z(x)
\ee
for some ${\bf C}_1^2$-valued function $z$.
Substituting (4.8) into the first equation of (4.6) gives 
\be
z^{\prime\prime}=-\mu\lambda z.
\ee
It follows from (4.4) that
\bea
\<z, z\>=\<z^{\prime}, z^{\prime}\>=0, \quad \<z, iz^{\prime}\>=\mu^{-1}.
 \eea
 Hence, $\phi$ is expressed as (4.1) satisfying (4.2).
Note that (4.9) is a consequence of (4.10).

Conversely, we can check that the immersion given in Theorem 7 is a quasi-biharmonic marginally trapped Lagrangian immersion
into ${\bf C}_1^2$.  The proof is finished.
\qed\vspace{1.5ex}

In the case when the ambient space is nonflat, we have
\begin{theorem}
Let $M$ be a marginally trapped Lagrangian surface
 in a $2$-dimensional
Lorentzian complex space form of constant holomorphic sectional curvature $4\epsilon\ne 0$.
Then, $M$ is quasi-biharmonic
if and only if the Gauss curvature of $M$ is equal to $\epsilon$.
\end{theorem}
\proof
Let $\phi: M\rightarrow \tilde M_1^2(4\epsilon)$ be a marginally trapped Lagrangian immersion
and let $\{e_1, e_2\}$ be  an orthonormal frame on $M$ satisfying
$\<e_1, e_1\>=1$, $\<e_2, e_2\>=-1$ and (3.17). Suppose that the second fundamental form is given by (3.2)-(3.3).

Assume that $\epsilon\ne 0$ and 
$M$ is quasi-biharmonic. Then,  it follows from (3.38)  that
$a-2b-c=0$ or
$(b-c)(a-2b-c)=3\epsilon$. 
If $(b-c)(a-2b-c)=3\epsilon$ holds, then
by using (3.19) we get
$G=2\epsilon$. 
On the other hand, by virtue of  Theorem 5.1 and 6.1 in \cite{chd}, we know that $G=0$ or $\epsilon$. 
Therefore, we obtain $\epsilon=0$, however this contradicts the assumption. Accordingly, we have
$a-2b-c=0$.  This and (3.19) show  $G=\epsilon$.

Conversely, assume that $G=\epsilon\ne 0$. Then, by (3.19) we have that $b=c$ or $a-2b-c=0$.
If $b=c$ holds,  by combining (3.13), (3.14), (3.18) and (3.24) we get 
\bea
\omega_1^2(e_1)=-\frac{e_1\alpha}{\alpha}, \quad \omega_1^2(e_2)=-\frac{e_2\alpha}{\alpha}.
\eea

On the other hand, the Gauss curvature $G$ is given by
\bea
G=-e_1(\omega_1^2(e_2))+e_2(\omega_1^2(e_1))+(\omega_1^2(e_1))^2-(\omega_1^2(e_2))^2.
\eea
Substituting (4.11) into (4.12), we have
\bea
G&=&\frac{1}{\alpha}\{e_1e_2\alpha-e_2e_1\alpha\}+\frac{1}{\alpha^2}\{(e_1\alpha)^2-(e_2\alpha)^2\}\nonumber\\
&=&\frac{1}{\alpha}[e_1, e_2]\alpha+\frac{1}{\alpha^2}\{(e_1\alpha)^2-(e_2\alpha)^2\}\nonumber\\
&=&\frac{1}{\alpha}\{\omega_1^2(e_1)e_1\alpha-\omega_1^2(e_2)e_2\alpha\}
+\frac{1}{\alpha^2}\{(e_1\alpha)^2-(e_2\alpha)^2\}\nonumber\\
&=&0, \nonumber
\eea
which contradicts the assumption. Hence, we have  $a-2b-c=0$.
Thus, it follows from (3.38) that $\tau_2(\phi)$ is lightlike, that is, $M$ is quasi-biharmonic
 because of $\epsilon\ne 0$. This completes the proof.
\qed\vspace{1.5ex}

By combining Theorem 5.1, 6.1 in \cite{chd} and Theorem 8, we obtain
\begin{theorem}
A quasi-biharmonic marginally trapped Lagrangian surface in ${\bf C}P_1^2(4)$ is 
locally congruent to the composition $\pi\circ L$, where $\pi: S^{5}_2(1)\rightarrow
{\bf C}P_1^2(4)$ is the Hopf fibration and $L$ is one of the following four families{\rm :}

{\rm (i)} 
\bea
L(x, y)=\frac{1}{a(x+y)}\Bigl(e^{\sqrt{2}iay}(\sqrt{2}+ia(x-y)), e^{\sqrt{2}iay}a(x-y), \sqrt{2}+ia(x+y)\Bigr),\nonumber
\eea
where $a$ is a nonzero real number.

{\rm (ii)} 
\bea
L(x, y)=\Bigl(\frac{2}{x+y}+\sqrt{2}if(y)\Bigr)z(y)-z^{\prime}(y),\nonumber
\eea
where $f(y)$ is a nonconstant real-valued function and $z(y)$ is a unit speed spacelike Legendre curve with 
the squared Legendre curvature $\hat\kappa^2=6f(y)^2$ in the light cone $\mathcal{L}$C.

{\rm (iii)} 
\bea
L(x, y)&=&\frac{e^{\frac{i}{\sqrt{2}}bx}}{3by}\Bigl((\sqrt{6}+i\sqrt{3}by)\cosh(\frac{\sqrt{6}}{2}bx)-
3(\sqrt{2}i+by)
\sinh(\frac{\sqrt{6}}{2}bx), \nonumber\\
&&3by\cosh(\frac{\sqrt{6}}{2}bx)
+\sqrt{3}(iby-2\sqrt{2})\sinh(\frac{\sqrt{6}}{2}bx),
\frac{\sqrt{6}+i\sqrt{3}by}{e^{3ibx/\sqrt{2}}}\Bigr),\nonumber
\eea
where $b$ is a positive real number.

{\rm (iv)} 
\bea
L(x, y)=\frac{z(y)}{x+y}-\frac{z^{\prime}(y)}{2},\nonumber
\eea
where $z(y)$ is a unit speed timelike special Legendre curve in the light cone $\mathcal{L}C$ with zero squared
Legendre curvature
and nonzero Legendre torsion.
\end{theorem}

\begin{theorem}
A quasi-biharmonic marginally trapped Lagrangian surface in ${\bf C}H_1^2(-4)$ is 
locally congruent to the composition $\pi\circ L$, where $\pi: H^{5}_3(-1)\rightarrow
{\bf C}H_1^2(-4)$ is the Hopf fibration and $L$ is one of the following four families{\rm :}

{\rm (i)} 
\bea
L(x, y)=\frac{1}{a(x-y)}\Bigl(ae^{\sqrt{2}iay}(x+y), a(x-y)+i\sqrt{2}, e^{\sqrt{2}iay}(a(x+y)+\sqrt{2}i)\Bigr),\nonumber
\eea
where $a$ is a nonzero real number.

{\rm (ii)} 
\bea
L(x, y)=\Bigl(\frac{2}{x-y}-\sqrt{2}if(y)\Bigr)z(y)+z^{\prime}(y),\nonumber
\eea
where $f(y)$ is a nonconstant real-valued function and $z(y)$ is a unit speed timelike Legendre curve with 
nonconstant squared Legendre curvature in the light cone $\mathcal{L}$C.

{\rm (iii)} 
\bea
&&L(x, y)=\frac{e^{\frac{-i}{\sqrt{2}}bx}}{3by}\Bigl(3(by+i\sqrt{2})\cosh(\frac{\sqrt{6}}{2}bx)+
3(\sqrt{2}+iby)
\sinh(\frac{\sqrt{6}}{2}bx), \nonumber\\
&&e^{\frac{3}{\sqrt{2}}ibx}(\sqrt{6}+i\sqrt{3}by), 3by\sinh(\frac{\sqrt{6}}{2}bx)
-\sqrt{3}(iby-2\sqrt{2})\cosh(\frac{\sqrt{6}}{2}bx)
\Bigr),\nonumber
\eea
where $b$ is a nonzero real number.

{\rm (iv)} 
\bea
L(x, y)=\frac{z(y)}{x-y}-\frac{z^{\prime}(y)}{2},\nonumber
\eea
where $z(y)$ is a unit speed spacelike special Legendre curve in the light cone $\mathcal{L}C$ with  
zero squared Legendre curvature
and nonzero Legendre torsion.
\end{theorem}


\section{Remarks {\normalsize(added on December 2, 2014)}}

\begin{remark}
{\rm By changing the sign of the metric of ${\bf C}P_1^2(4)$, we have that 
${\bf C}P_1^2(4)$ is holomorphically anti-isometric  to ${\bf C}H_1^2(-4)$.
Hence, it is sufficient to  consider only the case of ${\bf C}P_1^2(4)$.}
\end{remark}

\begin{remark}{\rm Comments on Theorem 5.1 in \cite{chd}:
Surfaces  in Case
 (A.b.ii) satisfy $f=0$. Hence, by (5.17) we have $\mu=0$.  This and (5.8) imply that 
these surfaces are minimal.
Accordingly,  (p.1.4) and (p.1.5) in  Theorem 5.1
should be removed from the list of marginally trapped surfaces.
Also, surfaces in Case (A.b.i)  are represented by (5.30) with (5.31), where
 $\delta=0$ and $f=a$ for some
nonzero constant $a$. 

Considering these facts, Theorem 9 of this paper should be read as follows.

\vskip10pt

\noindent {\bf Theorem 9}\hskip 5pt {\it A quasi-biharmonic marginally trapped Lagrangian surface in ${\bf C}P_1^2(4)$ is 
locally congruent to the composition $\pi\circ L$, where $\pi: S^{5}_2(1)\rightarrow
{\bf C}P_1^2(4)$ is the Hopf fibration and $L$ is given by 
\bea
L(x, y)=\Bigl(\frac{2}{x+y}+\sqrt{2}if(y)\Bigr)z(y)-z^{\prime}(y).\nonumber
\eea
Here $f(y)$ is a  real-valued function and $z(y)$ is a unit speed spacelike Legendre curve in the light cone $\mathcal{L}$C which satisfies
\bea
z^{\prime\prime\prime}=2\sqrt{2}ifz^{\prime\prime}+2(f^2+\sqrt{2}if^{\prime})z^{\prime}
+(\sqrt{2}i(f^{\prime\prime}+2\delta)+2ff^{\prime})z \nonumber
\eea
for some real-valued function $\delta(y)$.

}

\vskip10pt

The squared Legendre curvature and the Legendre torsion of 
$z$ described in Theorem 9 are given by $6f^2$ and 
$-8\sqrt{2}f^3+\sqrt{2}(f^{\prime\prime}+2\delta)$, respectively (see the proof of
Case (A.b.iii) in \cite{chd}).}

\end{remark}


\begin{remark}
{\rm In ["An isometric embedding of the complex hyperbolic space in a pseudo-Euclidean
space and its application to the study of real hypersurfaces", Tsukuba J. Math. {\bf 14}
(1990), 293-313],
Garay and Romero constructed an isometric embedding $\Psi$ of
${\bf C}H^n(-4)$ into some pseudo-Euclidean space, which is provided with 
good geometric properties as the first standard embedding of ${\bf C}P^n(4)$. 
 
 Let $\phi: M\rightarrow {\bf C}H^n(-4)$ be a real hypersurface with constant mean curvature.
 Then, they showed that 
 $\tau_2(\Psi\circ\phi)=Q$ holds for some non-zero constant vector if and only if 
 $M$ is locally congruent to 
 the horosphere $\pi(\{(z_1, z_2, \ldots, z_{n+1})\in H^{2n+1}_1(-1): |z_1-z_2|^2=1\})$. 
In this case, $Q$ is lightlike, i.e., the isometric
 immersion $\Psi\circ\phi$ is quasi-biharmonic.}
\end{remark}



\renewcommand{\baselinestretch}{1}

{\renewcommand{\baselinestretch}{1}

\bigskip

\noindent  General Education and Research Center\\
 Hachinohe Institute of Technology\\
Hachinohe 031-8501\\
       Japan\\
E-mail: sasahara@hi-tech.ac.jp


\begin{thebibliography}{99}

\bibitem{bal} A. Balmu\c{s}, {\it Biharmonic maps and submanifolds}, Balkan Society of Geometers, Differential Geometry-Dynamical Systems,
Monographs, 2009.




\bibitem{ba 2} M. Barros and A. Romero, {\it Indefinite Kaehler
manifolds,}
 Math. Ann. {\bf 261} (1982), 55-62.







\bibitem{ch1} B. Y. Chen, {\it Total Mean Curvature and Submanifold of
Finite Type,} World Scientific Publ., Singapore,  1984.



\bibitem{ch2} B. Y. Chen, {\it Maslovian Lagrangian surfaces of constant curvature
in complex projective or complex hyperbolic planes,} Math. Nachr. {\bf 278} (2005), 1242-1281.

\bibitem{ch3} B. Y. Chen, {\it Black holes, marginally trapped surfaces and quasi-minimal surfaces}
 Tamkang J. Math. {\bf 40} (2009), 313-341.






\bibitem{chd} B. Y. Chen and F. Dillen, {\it
Classification of marginally trapped Lagrangian surfaces
in Lorentzian complex space forms}, J. Math. Phys. {\bf 48} (2007), 013509, 23 pp;
Erratum, J. Math. Phys. {\bf 49} (2008), 059901.
















\bibitem{MO} S. Montaldo and C. Oniciuc {\it A short survey on biharmonic maps between Riemannian manifolds,}
Rev. Un. Mat. Argentina {\bf 47} (2006), no. 2, 1-22. 
Available at \texttt{http://inmabb.criba.edu.ar/revuma/}






\bibitem{sasa1} T. Sasahara,
{\it Quasi-minimal Lagrangian surfaces whose mean curvature vectors are eigenvectors.}
Demonstratio Math. {\bf 38} (2005), 185-196. 


\bibitem{sasa2}
T. Sasahara,
{\it Biharmonic Lagrangian surfaces of constant mean curvature in 
complex space forms,}
Glasgow Math. J. {\bf 49} (2007), 487-507.


















\end{thebibliography}
 \end{document}